# Exponential stability of second order delay differential equations through Floquet theory

Alexander Domoshnitsky*, Sergey Malev †, Tsahi Shavit ‡

**Abstract**

In this paper, we obtain results on exponential stability of second order delay differential equations, which are based on a version of the Floquet theory for delay differential equations of the second order we proposed. Our version allows researchers to preserve the order of equation and to obtain analogues of the classical results of the Floquet theory known for ordinary differential equations. On the basis of our version of the Floquet theory, new original unexpected results on the exponential stability are proposed. We demonstrate that choosing period of coefficients and delays of the gain in corresponding intervals allows to achieve the exponential stabilization in the cases considered as impossible when the standard technique was applied.

# Introduction

The Floquet theory is one of recognized areas in the qualitative theory of ordinary differential equations started with the classical paper by Floquet [23]. The Floquet theory provides formulas of solutions' representation. It makes this theory very convenient in stability studies. In this paper, we study stability of delay differential equations of the second order on the basis of the Floquet theory.

The use of equations with memory in applications leads to building analogues of the Floquet theory for equations with memory. It was developed in the following directions. The foundations of the Floquet theory for delay equations were proposed in the classical books [27, 30]. The foundations of the Floquet–Lyapunov theory for discrete time systems with periodic coefficients are presented in [28, 29, 36]. Lyapunov-Floquet representations for differential equations with delay were proposed in [39]. Note several papers on applications of the Floquet theory to delay equations: an application to Bloch's theorem for nonlocal potentials was proposed in [50], a cholera epidemic model with phage dynamics and seasonality was considered in [49], an application to the control of drone's flight see in [11]. The Floquet theory for equations with distributed delays was proposed in [2, 3, 9, 10], where results on stability of these equations were obtained.

---

*Ariel University, Ariel, Israel, adom@ariel.ac.il
†Ariel University, Ariel, Israel,
‡Ariel University, Ariel, Israel,

The classical definition of homogeneous equation [42, 30] and its solutions as continuation of every possible initial functions leads us to infinite-dimensional fundamental system. It should be stressed that an infinite-dimensional fundamental system does not allow obtaining a full analogue of the Floquet theory as it was noted in the well known book [31] (see p. 236). In this connection, we can note, simple solutions' representations and assertions about passing to systems with constant coefficients which are so important in the Floquet theory for ordinary differential equations (ODEs).

We propose a version of the Floquet theory for delay differential equations. In our version, the orders of equation and the monodromy matrix coincide. This makes our version similar to the classical Floquet theory for ODEs [23] and allows us to use the classical methods of the theory developed for ODEs.

A principal development of our approach is to preserve not only finite-dimensional fundamental systems of ordinary delay differential equations, but also the number of independent solutions in it. The idea to come to finite-dimensional fundamental system for delay equations and to propose an analog of the Floquet theory on the idea of "finite-dimensional basis" was developed in the recent paper [47]. Initial functions are presented there in the form of linear combination of corresponding functions. The Floquet solutions are constructed as continuations of these functions. Semi discretization concept, where sampled data approach was combined with periodicity was presented in [32] and a way to representation of solution was obtained. An idea to come to finite-dimensional monodromy operator was presented in [34]. Although the fundamental system in the works [32, 34, 33, 47] is finite-dimensional, the Floquet solutions are presented as a solutions of a corresponding non-homogeneous equation, and as a result, the dimension of the monodromy operator is greater than the order of the given equation. This implies that the properties of solutions of equations that are essentially related to its order are lost. The approach, known as act-and-wait concept by Stepan and Insperger, where the orders of equation and monodromy operators coincide, was described in [48]. Note the paper [52], where the idea of act-and-wait control was developed for delay systems. Our approach can be considered, in corresponding sense, as development of act-and-wait control approach.

Another idea is stabilizing effect of delays when equations without delays are unstable but adding delays becomes them to stable. The idea of stabilization by appearing delays in corresponding terms of unstable equations was noted as important in [44], where corresponding references can be found. Note in this connection the paper [15], where it was demonstrated that although the ordinary differential equation

$$x''(t) + \left\{\sum_{i=1}^{m} b_i(t)\right\} x(t) = 0, \ t \in [0, \infty)$$

can be oscillatiting and asymptoticaly unstable, the delay equation

$$x''(t) + \sum_{i=1}^{m} b_i(t) x(t - \tau_i(t)) = 0, \ t \in [0, \infty) \tag{0.1}$$

under corresponding conditions on coefficients $b_i(t)$ and delays $\tau_i(t)$ is nonoscillating and exponentially stable. The basic idea in [15] is to avoid the condition on nonnegativity of the coefficients $b_i(t)$ for all $i = 1, ..., m$, and to allow terms with positive and terms with negative coefficients $b_i(t)$ to compensate each other. Note that one-term equation with constant coefficient and delay

$$x''(t) + bx(t-\tau) = 0, \ t \in [0, \infty)$$

has unbounded solution for every positive $b$ and $\tau$ [42, 30]. For equation with a variable delay $\tau(t)$, the condition

$$\int_0^\infty \tau(t)dt < \infty$$

is necessary and sufficient for boundedness of all solutions on the semi axis $t \geq 0$ [14]. In the case of advanced argument, the conditions of boundedness of solutions were obtained in [22]. In the case of two-term autonomous equation ($m = 2$) with $\tau_1 = 0$ and $b_1 b_2 < 0$ conditions of the exponential stability were obtained in [13]. For equation with damping term (1.1) conditions of exponential stability in the case of autonomous equation were can be found in [13], for equations with constant delay in the term with derivative or without delay in derivative were obtained in [37, 38], for variable delay in the term with derivative, the first result on exponential stability was obtained in [18] and various results and a concept of the stability studies were proposed in [8]. Other results on stability of second order equations can be found in [13, 15, 22].

In our paper, we use also the idea of sampled-data controllers [24, 25, 32], but in a direction to come to a version of the Floquet theory. A goal of our approach in building a new version of the Floquet theory for delay equations is to remain with a finite-dimensional fundamental systems of homogeneous equations.

We base our approach on Azbelev's definition of the homogeneous equations. This allows us to stay with finite-dimensional fundamental system [4, 5] dimension of which is equal to the order of the equation.

Our paper is built as follows. In Section 1 we explain what is our direction in developement of act-and-wait control. In Section 2, some facts of the general theory of delay equations (see [4, 5]) are noted. The concept of the Floquet theory was proposed in Section 3. Various examples of unexpected results on the uniform exponential stability were proposed in Section 4. Explanations of examples are proposed in Section 5. We build programs based on described algorithms of the exponential stability based on the Floquet theory and act-and-wait concept formulated in Section 3. Discussion on conclusions and open problems can be found in Section 6.

# 1 Developement of the idea of act-and-wait control

In this paper, we consider the equation

$$x''(t) + a(t)x'(t) + \sum_{i=1}^m b_i(t)x(t-\tau_i(t)) = f(t), \ t \in [0, \infty) \tag{1.1}$$

$$x(\xi) = \varphi(\xi) \text{ for } \xi < 0, \tag{1.2}$$

with essentially bounded Borel measurable right-hand side $f(t)$, initial function $\varphi(t)$, coefficients $a(t), b_i(t)$ and delays $\tau_i(t) \geq 0$ for $t \in [0, \infty)$.

How could we define its homogeneous equation? The traditional approach [42, 30] defines the homogeneous equation as

$$x''(t) + a(t)x'(t) + \sum_{i=1}^{m} b_i(t)x(t - \tau_i(t)) = 0, \ t \in [0, \infty) \qquad (1.3)$$

with the initial functions defined by (1.2). If a researcher wants to consider solutions as continuation of every possible initial function $\varphi(t)$, then the infinite-dimensional fundamental system appears as a result of this definition of the homogeneous equation. The idea of N.V.Azbelev [4, 5] was to define the homogeneous equation as (1.3) with the zero initial function

$$x(\xi) = 0 \text{ for } \xi < 0. \qquad (1.4)$$

It should be noted that this zero initial function allows to consider all solutions of equation (1.1),(1.2) for all possible initial functions $\varphi(t)$. Actually, we can rewrite equation (1.1) in the form

$$\begin{aligned} x'(t) + a(t)x(t) + \sum_{i=1}^{m} b_i(t)\sigma(t - \tau_i(t))x(t - \tau_i(t)) = \\ f(t) - \sum_{i=1}^{m} b_i(t)\varphi(t - \tau_i(t))\left\{1 - \sigma(t - \tau_i(t))\right\}, \end{aligned} \quad t \in [0, \infty) \qquad (1.5)$$

where

$$\sigma(t - \tau_i(t)) = \begin{cases} 1, & t - \tau_i(t) \geq 0 \\ 0, & t - \tau_i(t) < 0 \end{cases}. \qquad (1.6)$$

The homogeneous equation for (1.1) can be considered in the form

$$x''(t) + a(t)x'(t) + \sum_{i=1}^{m} b_i(t)\sigma(t - \tau_i(t))x(t - \tau_i(t)) = 0, \ t \in [0, \infty) \qquad (1.7)$$

or we can write it as

$$x''(t) + a(t)x'(t) + \sum_{i=1}^{m} b_i(t)x(t - \tau_i(t)) = 0, \ t \in [0, \infty) \qquad (1.8)$$

with zero initial function (1.4).

**Remark 1.1.** Stability of (1.1) with respect to right-hand side, which is most important type of stability in applications, is equivalent (in the case of bounded delays), according to Bohl-Perron theorem (see, for example, pp. 499-501 of [1]), to the uniform exponential stability (see Definition 2.1 below) of homogeneous equation (1.7) or of equation (1.8) with initial function (1.4).

The concept of finite-dimensional fundamental system allows N.V.Azbelev and his followers to build a theory of functional differential equations presented in the books [4, 6, 5, 1, 8]. Stability is one of the main parts of this theory. All results on stability obtained in its frame are based on this definition of homogeneous equations.

It is surprising that the same idea of a finite-dimensional fundamental system for equations with delays came from sufficiently practical ideas in the field of control. As very important example, the concept of so-called act-and-wait control [48] can be noted. Let us try to describe one of its basic ideas. The following simple equation was considered

$$x''(t) + ax'(t) + B(t)x(t - \tau) = 0, \ t \in [0, \infty), \tag{1.9}$$

where

$$B(t) = \begin{cases} 0, & 0 \le t < \tau \\ b, & \tau \le t < \omega \end{cases} ; \ B(t) = B(t - \omega) \text{ for } t \in [\omega, \infty), \tag{1.10}$$

with constant $a, b$ and $\tau$. The control term $B(t)x(t - \tau)$ is active only on the intervals $[\tau, \omega], [\tau + \omega, 2\omega], \ldots, [\tau + n\omega, (n+1)\omega], \ldots$. They are "act" intervals, others - "wait" ones. From here comes the name of the method: act-and-wait control. The goal is to avoid dependence of the solution on the initial function in order to stay with two-dimensional fundamental system. The "wait" intervals can be problematic with the point of view of many applications (see, for example, projects on vision-based drone navigator [16] and on navigation of ground robot from upper position [17] based on the models described in [7, 11]). Is it possible to arrive at a finite-dimensional equation without assumption on existence of "wait" intervals in the control $B(t)x(t - \tau)$? The idea of E.Fridman and her co-authors [25, 24] to use the delay in the form $\tau(t) = t$, for $t \in [0, \omega)$, $\tau(\omega) = 0$ and then to be continued as $\omega$−periodic function, can be considered as a way to achieve this. This idea open another way to come to finite-dimensional delay equations. A generalized form of this idea was obtained independenltly and presented as a concept in [11, 19, 21] for first order delay equations.

Thus, one of the possible ways to come to finite-dimensional fundamental system is to change not the coefficients in (1.3), but the delays. In equation (1.8), we can assume that

$$t - \tau_i(t) \ge 0, \ t \in [0, \omega], \ i = 1, ..., m. \tag{1.12}$$

There is no now any influence of the initial function on the solution $x(t)$ of equation (1.8). Thus, we come to two-dimensional delay equation (1.8). To obtain a version of the Floquet theory, we assume periodic continuation of the coefficients and delays

$$a(t) = a(t - \omega), \ b_i(t) = b_i(t - \omega), \ \tau_i(t) = \tau_i(t - \omega) \text{ for } t \in [\omega, \infty). \tag{1.13}$$

We can continue from this point and not from the scheme of "act-and-wait" control only. Consider, for example, the equation

$$x''(t) + ax'(t) + bx(t - \tau(t)) = 0, \ t \in [0, \infty) \tag{1.14}$$

where

$$\tau(t) = \begin{cases} t, & 0 \le t < \tau \\ \tau, & \tau \le t < \omega \end{cases} ; \ \tau(t) = \tau(t - \omega) \text{ for } t \in [\tau, +\infty), \tag{1.15}$$

where $a, b$ and $\tau$ are constants. Equation (1.14) looks similar to (1.9) and both of them look similar to the autonomous equation with constant coefficients and delay

$$x''(t) + ax'(t) + bx(t - \tau) = 0, \ t \in (-\infty, +\infty), \tag{1.16}$$

which is considered as a basic in the theory of delay differential equations. But equations (1.9) and (1.14) are one-dimensional, and, as a result, their asymptotic properties and, of course, the conditions of stability can be radically different from ones of infinite-dimensional equation (1.16). Studying equations (1.9) and (1.14) for coefficients $a, b$ and $\tau$ in the zones of instability (obtained, for example, through D-subdivision method [37]) of equation (1.16), we see that equations with small changes in delays can be exponentially stable. Holding $\omega$ as a parameter, allowed us to see a general qualitative picture. We have taken some steps towards understanding the qualitative picture to achieve stabilization. We obtain that in many cases there are infinite number of intervals $[\alpha_n, \beta_n]$, where $\alpha_n \to \infty$ when $t \to \infty$, such that if $\omega$ is situated in one of them, equations (1.9) or (1.14) are exponentially stable. The length of these intervals may tend to zero. If we don't know for sure from pure theoretical ideas that these intervals exist, it is almost impossible to find this property.

In our approach the order of equation and the number of solutions in its fundamental system coincide. For second order equations they are to two-dimensional. This is, in our opinion, an important advantage of our approach allowing to preserve specific properties of second order equations. Our simple version of the Floquet theory is completely analogous to the case of ODEs.

On the basis of our version of the Floquet theory, new original unexpected results on the exponential stability are proposed. We demonstrate that choosing period of coefficients and delays of the gain in corresponding intervals allows to achieve the exponential stabilization in the cases considered as impossible when the standard technique was applied. Our goal is to propose new possibilities of stabilization in difficult for the standard stability studies cases on the basis of the Floquet theory. In the frame of our approach, only the value of solutions and its first derivatives at the point $\omega$ (which is equal to the period of coefficients and delays) turns out to be important in stability analysis. We hope, that our approach develops ideas of stabilizing effects of delays [44] and add new possibilities for stability studies developed in the works [32, 34, 48, 47].

## 2 Preliminaries

Let us consider the equation

$$x''(t) + a(t)x'(t) + \sum_{i=1}^{m} b_i(t)x(t - \tau_i(t)) = f(t), \ t \in [0, \infty) \tag{2.1}$$

$$x(\xi) = 0, \ \xi < 0, \tag{2.2}$$

with essentially bounded measurable coefficients $b_i(t)$ and delays $\tau_i(t) \geq 0$ for $t \in [0, \infty)$. Its general solution was obtained in the paper [40] (then it also appeared in the books [4, 5])

$$x(t) = \int_0^t C(t, s)f(s)ds + x_1(t)x(0) + x_2(t)x'(0), \ t \in [0, \infty), \tag{2.3}$$

where $x_1(t)$ and $x_2(t)$ are two independent solutions of the homogeneous equation

$$x''(t) + a(t)x'(t) + \sum_{i=1}^{m} b_i(t)x(t - \tau_i(t)) = 0, \ t \in [0, \infty) \tag{2.4}$$

with the zero initial function (1.4) such that

$$x_1(0) = 1, \ x_1'(0) = 0, \ x_2(0) = 0, \ x_2'(0) = 1, \tag{2.5}$$

We use below the Wronskian of the fundamental system

$$W(t) = \begin{vmatrix} x_1(t) & x_2(t) \\ x_1'(t) & x_2'(t) \end{vmatrix}. \tag{2.6}$$

The kernel $C(t,s)$ of the integral term in this representation is called the Cauchy function (fundamental function in another termminology) of equation (2.1),(2.2). It is defined as the solution of so-called "$s$-truncated" homogeneous equations

$$x''(t) + a(t)x'(t) + \sum_{i=1}^{m} b_i(t)x(t - \tau_i(t)) = 0, \ t \in [s, \infty), \tag{2.7}$$

$$x(\xi) = 0, \ \xi < s, \tag{2.8}$$

for every fixed $s \geq 0$, satisfying the initial condition $x(s) = 0, x'(s) = 1$.

**Definition 2.1.** *We will say that equation (2.4) is uniformly exponentially stable if there exist $N > 0$ and $\gamma > 0$ such that the solution of problem (2.4) with the initial function*

$$x(\xi) = \psi(\xi) \text{ for } \xi < t_0, \tag{2.9}$$

*where $\psi : (-\infty, t_0)$ is Borel measurable bounded function, satisfies the estimate*

$$|x(t)| \leq Ne^{-\alpha(t-t_0)} \mathrm{esssup}_{\xi \leq t_0} |\psi(\xi)|, \ t \geq t_0, \tag{2.10}$$

*where the constants $N$ and $\alpha$ do not depend on $t_0$ and $\psi$.*

**Definition 2.2.** *We will say that the Cauchy function of equation (2.1) satisfies the exponential estimate if there exist such positive $N$ and $\alpha$ that*

$$|C(t,s)| \leq Ne^{-\alpha(t-s)} \text{ for } 0 \leq s \leq t < \infty. \tag{2.11}$$

**Remark 2.1.** It is known that for every bounded delays these two definitions are equivalent [4, 6].

**Remark 2.2.** It is clear that the stability with respect to right-hand side follows from the exponential estimate (2.11) of the Cauchy function $C(t,s)$.

# 3 Method: Floquet Theory for Delay Differential Equation of the Second Order

In this section, we propose a version of the Floquet theory for delay differential equations of the second order. The coefficients and delays in the frame of the Floquet theory are assumed to be periodic with a corresponding period $\omega$, i.e., the following conditions are true.

**Assumption 1:**

$$a(t + \omega) = a(t), \ b_i(t + \omega) = b_i(t), \ \tau_i(t + \omega) = \tau_i(t), \ t \in [0, \infty), \tag{3.1}$$

To preserve a finite-dimensional version of the Floquet theory we will need an additional conditions on the delays.

**Assumption 2:**
$$t - \tau_i(t) \geq 0, \ t \in [0, \omega], \ i = 1, ..., m. \tag{3.2}$$

In the frame of the Floquet theory we try to obtain existence of solution satisfying the equality
$$x(t + \omega) = \lambda x(t). \tag{3.3}$$

**Theorem 3.1.** *Let Assumptions 1 and 2 be fulfilled. Then solutions of equation (2.4) satisfying equality (3.3) exist and the numbers $\lambda$ are the roots of the following quadratic equation*
$$\lambda^2 - (x_1(\omega) + x_2'(\omega))\lambda + W(\omega) = 0. \tag{3.4}$$

This assertion can be considered as an analog of the Floquet theorem for delay differential equation (2.4).

**Proof.** Let us consider the set $Z$ of functions $z : [0, \infty) \to \mathbb{R}$ with absolutely continuous derivative $z'$ on every finite interval $[0, T]$ such that every function $z \in Z$ can be presented in the form
$$z(t + k\omega) = c_{1k}x_1(t) + c_{2k}x_2(t), \ k = 0, 1, 2, \ldots, \ t \in [0, \omega], \tag{3.5}$$

where $x_1(t), x_2(t)$ are two independent solution of (2.4) satisfying the initial conditions (2.5) and $c_{1k}, c_{2k}$ are constants. Assumptions 1 and 2 imply the existence of such solutions that (3.3) is fulfilled, according to [45] or [46], pp. 268-286.

Now let us obtain equation (3.4) for finding $\lambda$. Considering equality (3.5) for $k = 1$, using existence of such solutions $x$, we come to the system
$$\begin{cases} x(t + \omega) = c_1 x_1(t) + c_2 x_2(t) \\ x'(t + \omega) = c_1 x_1'(t) + c_2 x_2'(t) \end{cases}. \tag{3.6}$$

Setting $t = 0$, and substituting into equality (3.3), we come to the system
$$\begin{cases} x(\omega) = \lambda(c_1 x_1(0) + c_2 x_2(0)) \\ x'(\omega) = \lambda(c_1 x_1'(0) + c_2 x_2'(0)) \end{cases}. \tag{3.7}$$

Now, substituting the initial conditions (2.5), and using (3.6) for $t = 0$, we come to the system
$$\begin{cases} \lambda c_1 = c_1 x_1(\omega) + c_2 x_2(\omega) \\ \lambda c_2 = c_1 x_1'(\omega) + c_2 x_2'(\omega) \end{cases}. \tag{3.8}$$

and then to the system
$$\begin{cases} c_1(x_1(\omega) - \lambda) + c_2 x_2(\omega) = 0 \\ c_1 x_1'(\omega) + c_2(x_2'(\omega) - \lambda) = 0 \end{cases}. \tag{3.9}$$

Its determinant is
$$\Delta(\omega) = \begin{vmatrix} x_1(\omega) - \lambda & x_2(\omega) \\ x_1'(\omega) & x_2'(\omega) - \lambda \end{vmatrix}. \tag{3.10}$$

The equality $\Delta(\omega) = 0$ leads us to quadratic equation (3.4).

The proof of Theorem 3.1 has completed.

**Remark 3.1.** If the roots $\lambda_1$ and $\lambda_2$ of equation (3.4) satisfy the inequalities $|\lambda_{1,2}| < 1$, then the formulas of solutions' presentation obtained, for example, in the paper [39] imply the uniform exponential stability of equation (2.4). Note that this fact can be also obtained directly with short explanations from the classical books [45], [46].

# 4 Examples

**Example 4.1.** Consider the equation

$$x''(t) + ax(t) + bx(t - \tau(t)) = 0, \ t \in [0, \infty) \qquad (4.1)$$

for: $a = 1$, $b = -0.5$, and

$$\tau(t) = \begin{cases} t, & 0 \leq t < 0.5 \\ 0.5, & 0.5 \leq t < \omega \end{cases} ; \ \tau(t) = \tau(t - \omega) \text{ for } t \in [\tau, +\infty). \qquad (4.2)$$

This equation is uniformly exponentially stable when $\omega$ is situated in one of the intervals $[0.51, 0.6]$ or $[2.8, 20.0]$. Note that [15] and Theorem 6.1 of [8] (see also Example 6.1 there on p.115) claim that the inequality $0 < 4(a - |b|) \leq \frac{1}{e^2(\max \tau(t))^2}$ implies the exponential stability of (4.1). In our case, this inequality is not fulfilled.

**Example 4.2.** Consider equation (4.1) for $a = 1, b = -1$ and $\tau(t)$ defined by (4.2).

This equation is uniformly exponentially stable for $\omega \in [0.51, 1.2]$. Note that in contrast with assertions in Chapter 6 of [8], the inequality $a > |b|$ is not fulfilled here.

**Example 4.3.** Consider the equation

$$x''(t) + ax'(t) + bx(t - \tau(t)) = 0, \ t \in [0, \infty) \qquad (4.3)$$

for: $a = 0.5, b = 1$ and $\tau(t)$ defined by (4.2). This equation is uniformly exponentially stable for $\omega \in [0.5, 6]$. Note that the uniform exponential stability of (4.3) was obtained in [18] through positivity-based approach. The main condition there was existence of positive $\alpha$ such that the inequality $\alpha^2 + be^{\alpha \max \tau(t)} \leq \alpha/2$ was fulfilled. Its conditions cannot be fulfilled since we come to the inequality $\alpha^2 - \frac{1}{2}\alpha + e^{\alpha/2} \leq 0$. Even in the case of "very small" $\alpha$, we cannot obtain that this inequality is fulfilled.

**Example 4.4.** Consider equation (4.3) with the delay defined as

$$\tau(t) = \begin{cases} t, & 0 \leq t < 1 \\ 1, & 1 \leq t < \omega \end{cases} ; \ \tau(t) = \tau(t - \omega) \text{ for } t \in [\tau, +\infty). \qquad (4.5)$$

In the case of this $\tau(t)$, equation (4.3) is uniformly exponential stable for
1) $a = 2, b = 2.5$, if the period $\omega$ is situated in one of the intervals $[2.05, 3.4], [5.0, 6.25]$, $[8.0, 9.05], [10.95, 11.9], [13.9, 14.75], [16.85, 17.6], [19.8, 20.4], [22.75, 23.25], [25.7, 26.1]$, $[28.65, 28.95], [31.65, 31.8]$;
2) $a = 2, b = 3$ if the period $\omega$ is situated in $[2.35, 2.75]$;
3) $a = 2, b = 3.3$, if the period $\omega$ is situated in $[2.45, 2.55]$;
4) $a = 2.5, b = 4.6$, if the period $\omega$ is situated in $[2.25, 2.3]$.
As well as we know, results [37, 38, 13, 8] does not work here.

# 5 Algorithms for Solving Equations: Explanations to Examples

The idea of our approach is to construct solutions $x_1(t), x_2(t)$ of delay equation (1.1) on corresponding interval and then to built equation (3.4). In the case of roots $\lambda_1$ and $\lambda_2$ satisfying inequalities $|\lambda_{1,2}| < 1$, we obtain the exponential stability of delay equations. Construction of solutions can be done as follows.

## 5.1 The equation $x''(t) + ax'(t) + bx(t - \tau(t)) = 0$.

In this section we will consider the equation $x''(t) + ax'(t) + bx(t - \tau(t)) = 0$, where

$$\tau(t) = \begin{cases} t, & 0 \leq t < \tau, \\ \tau, & t \geq \tau. \end{cases} \quad (5.1)$$

The function $x(t)$ is defined for $t \geq 0$, with two different initial conditions, namely,

$$\begin{cases} x(0) = 0 \\ x'(0) = 1 \end{cases} \text{ and } \begin{cases} x(0) = 1 \\ x'(0) = 0 \end{cases}.$$

We will consider functions $X_n(t) : [0, \tau] \to \mathbb{R}$ defined on the interval $[0, \tau]$ which define our function $x(t)$ in the spiral way, that is

$$x(t) = X_n(t - n\tau), \ t \in [n\tau, (n+1)\tau], \ n \geq 0.$$

For $n = 0$ we have $X_0''(t) + aX_0'(t) = -bX_0(0)$ with two different initial conditions. The solution is $X_0(t) = A + Bt + Ce^{-at}$ where in the first task $A = \frac{1}{a}, B = 0, C = -\frac{1}{a}$ and in the second task $A = 1 + \frac{b}{a^2}, B = -\frac{b}{a}, C = -\frac{b}{a^2}$.

For $n > 0$ the equation is $X_n''(t) + aX_n'(t) = -bX_{n-1}(t)$ with two initial conditions: $X_n(0) = X_{n-1}(\tau)$ and $X_n'(0) = X_{n-1}'(\tau)$.

The solution will have a form $X_n(t) = P_n(t) + Q_n(t)e^{-at}$, where $P_n(t)$ is a polynomial of degree no more than $n + 1$ and $Q_n(t)$ is a polynomial of degree no more than $n$. We can define

$$P_n(t) = \sum_{i=0}^{n+1} \alpha[n, i]t^i, \ Q_n(t) = \sum_{i=0}^{n} \beta[n, i]t^i.$$

Note, $X_n'(t) = P_n'(t) + (Q_n'(t) - aQ_n(t))e^{-at}$ and $X_n''(t) = P_n''(t) + (Q_n''(t) - 2aQ_n'(t) + a^2Q_n(t))e^{-at}$. Hence, $X_n''(t) + aX_n'(t) = P_n''(t) + aP_n'(t) + (Q_n''(t) - aQ_n'(t))e^{-at}$, which equals, as we know, to $-bX_{n-1}(t) = -bP_{n-1}(t) - bQ_{n-1}(t)e^{-at}$. Thus, $P_n''(t) + aP_n'(t) = -bP_{n-1}(t)$ and $Q_n''(t) - aQ_n'(t) = -bQ_{n-1}(t)$. Remind,

$$P_n(t) = \sum_{i=0}^{n+1} \alpha[n, i]t^i,$$

hence

$$P_n'(t) = \sum_{i=0}^{n} (i+1)\alpha[n, i+1]t^i,$$

and

$$P_n''(t) = \sum_{i=0}^{n-1} (i+2)(i+1)\alpha[n, i+2]t^i.$$

Therefore,

$$P_n''(t) + aPn'(t) = \sum_{i=0}^{n-1} ((i+2)(i+1)\alpha[n, i+2] + a(i+1)\alpha[n, i+1])t^i + a(n+1)\alpha[n, n+1]t^n,$$

and it equals to
$$-bP_{n-1}(t) = \sum_{i=0}^{n} -b\alpha[n-1, i]t^i.$$

Therefore, $\alpha[n, n+1] = -\frac{b\alpha[n-1,n]}{a(n+1)}$. Now let $k$ run from $n$ down to 1, we have
$$(k+1)k\alpha[n, k+1] + ak\alpha[n, k] = -b\alpha[n-1, k-1]$$
and thus we compute
$$\alpha[n, k] = \frac{1}{ak}\left(-b\alpha[n-1, k-1] - (k+1)k\alpha[n, k+1]\right)$$

The coefficients of $Q_n(t)$ we compute using $Q_n''(t) - aQ_n'(t) = -bQ_{n-1}(t)$.

The polynomial $Q_n(t)$ equals $\sum_{i=0}^{n} \beta[n, i]t^i$, hence
$$Q_n'(t) = \sum_{i=0}^{n-1} (i+1)\beta[n, i+1]t^i,$$
and
$$Q_n''(t) = \sum_{i=0}^{n-2} (i+2)(i+1)\beta[n, i+2]t^i.$$

Hence,
$$Q_n''(t) - aQ_n'(t) = \sum_{i=0}^{n-2} \left((i+2)(i+1)\beta[n, i+2] - a(i+1)\beta[n, i+1]\right)t^i + (-an)\beta[n,n]t^{n-1}$$
which equals to
$$-bQ_{n-1}(t) = \sum_{i=0}^{n-1} -b\beta[n-1, i]t^i.$$

Therefore, $\beta[n, n] = \frac{b\beta[n-1,n-1]}{an}$. Now let $k$ run from $n-1$ down to 1, we have $(k+1)k\beta[n, k+1] - ak\beta[n, k] = -b\beta[n-1, k-1]$, hence
$$\beta[n, k] = \frac{b\beta[n-1, k-1] + (k+1)k\beta[n, k+1]}{ak}.$$

Now we will use the initial conditions $X_n(0) = X_{n-1}(\tau)$ and $X_n'(0) = X_{n-1}'(\tau)$. As soon as $X_n(t) = P_n(t) + Q_n(t)e^{-at}$,
$$X_n(0) = P_n(0) + Q_n(0) = \alpha[n, 0] + \beta[n, 0] = X_{n-1}(\tau) \tag{5.2}$$
as well as
$$X_n'(0) = P_n'(0) + Q_n'(0) - aQ_n(0) = X_{n-1}'(\tau). \tag{5.3}$$

The equation (5.3) gives
$$\beta[n, 0] = Q_n(0) = \frac{-X_{n-1}'(\tau) + P_n'(0) + Q_n'(0)}{a}.$$

And the equation (5.2) gives
$$\alpha[n, 0] = P_n(0) = X_{n-1}(\tau) - \beta[n, 0].$$

Thus, the recurrence formula has been obtained.

## 5.2 Equation $x''(t) + ax(t) + bx(t - \tau(t)) = 0.$

In this section, we shall solve the equation

$$x''(t) + ax(t) + bx(t - \tau(t)) = 0, \quad (a > 0)$$

where

$$\tau(t) = \begin{cases} t, & 0 \leq t < \tau, \\ \tau, & t \geq \tau. \end{cases} \tag{5.4}$$

with two different boundary conditions, namely

$$\begin{cases} x(0) = 0 \\ x'(0) = 1 \end{cases} \text{ and } \begin{cases} x(0) = 1 \\ x'(0) = 0 \end{cases}.$$

Here we again make use of the spiral representation of the function

$$X(t) = X_n(t - n\tau), \quad t \in [n\tau, (n+1)\tau], \ n \geq 0.$$

As in the previous cases, all functions $X_n(t) : [0, \tau] \to \mathbb{R}$ are defined on the interval $[0, \tau]$, and the dependence is obtained through the equation

$$X_n''(t) + aX_n(t) = -bX_{n-1}(t).$$

That is, at each step we must solve a harmonic differential equation.

At the initial stage, we obtain the equation

$$X_0''(t) + aX_0(t) = -bX_0(0)$$

with the corresponding initial conditions. In the first problem we obtain

$$\begin{cases} X_0''(t) + aX_0(t) = 0 \\ X(0) = 0 \\ X'(0) = 1 \end{cases},$$

whose solution is the function $X_0(t) = \frac{1}{\sqrt{a}} \sin(\sqrt{a}t)$.

In the second problem we obtain

$$\begin{cases} X_0''(t) + aX_0(t) = -b \\ X(0) = 1 \\ X'(0) = 0 \end{cases},$$

whose solution is $X_0(t) = \left(\frac{b}{a} + 1\right) \cos(\sqrt{a}t) - \frac{b}{a}$.

Thus, both problems can be unified into one by assuming that

$$X_0(t) = A\cos(\sqrt{a}t) + B\sin(\sqrt{a}t) + C,$$

taking for the first problem $A = C = 0, B = \frac{1}{\sqrt{a}}$, and for the second $A = \frac{b}{a} + 1, B = 0, C = -\frac{b}{a}$.

Hence, the problem is as follows: given the function

$$X_0(t) = A\cos(\sqrt{a}t) + B\sin(\sqrt{a}t) + C,$$

for a triple of coefficients $A, B, C$, we must solve the sequence of differential equations

$$\begin{cases} X_n''(t) + aX_n(t) = -bX_{n-1}(t) \\ X_n(0) = X_{n-1}(\tau) \\ X_n'(0) = X_{n-1}'(\tau) \end{cases}.$$

Such equations frequently arise in the general course on differential equations, and it is known that the solution has the form

$$X_n(t) = P_n(t)\cos(\sqrt{a}t) + Q_n(t)\sin(\sqrt{a}t) + C_n,$$

where $P_n, Q_n$ are polynomials of degree at most $n$, and $C_n$ is a constant. The polynomials are defined as

$$P_n(t) = \sum_{i=0}^{n} \alpha[n,i]t^i, \quad Q_n(t) = \sum_{i=0}^{n} \beta[n,i]t^i.$$

Since $X_n(t) = P_n(t)\cos(\sqrt{a}t) + Q_n(t)\sin(\sqrt{a}t) + C_n$,

$$X_n'(t) = (P_n'(t) + Q_n(t)\sqrt{a})\cos(\sqrt{a}t) + (Q_n'(t) - P_n(t)\sqrt{a})\sin(\sqrt{a}t),$$

and therefore

$$X_n''(t) = (P_n''(t) + 2Q_n'(t)\sqrt{a} - aP_n(t))\cos(\sqrt{a}t) + (Q_n''(t) - 2P_n'(t)\sqrt{a} - aQ_n(t))\sin(\sqrt{a}t).$$

Consequently,

$$X_n''(t) + aX_n(t) = (P_n''(t) + 2Q_n'(t)\sqrt{a})\cos(\sqrt{a}t) + (Q_n''(t) - 2P_n'(t)\sqrt{a})\sin(\sqrt{a}t) + aC_n.$$

From the equation it follows that this equals

$$-bX_{n-1}(t) = -bP_{n-1}(t)\cos(\sqrt{a}t) - bQ_{n-1}(t)\sin(\sqrt{a}t) - bC_{n-1}.$$

Hence it follows immediately that

$$-\tfrac{b}{a}C_{n-1} = C_n,$$

that is, $C_n = \left(-\tfrac{b}{a}\right)^n C$, which requires no recurrence relation.

Furthermore, we have

$$\begin{cases} P_n''(t) + 2\sqrt{a}Q_n'(t) = -bP_{n-1}(t) \\ Q_n''(t) - 2\sqrt{a}P_n'(t) = -bQ_{n-1}(t) \end{cases}.$$

This is a system of differential equations with respect to the polynomials $P_n'(t)$ and $Q_n'(t)$. For convenience, introduce the notation $R(t) = P_n'(t)$ and $S(t) = Q_n'(t)$. Then

$$\begin{cases} R'(t) + 2\sqrt{a}S(t) = -bP_{n-1}(t) \\ S'(t) - 2\sqrt{a}R(t) = -bQ_{n-1}(t) \end{cases}.$$

From the first line we obtain

$$S(t) = \tfrac{1}{2\sqrt{a}}(-bP_{n-1}(t) - R'(t)).$$

Substituting into the second gives

$$\tfrac{1}{2\sqrt{a}}(-bP'_{n-1}(t) - R''(t)) - 2\sqrt{a}R(t) = -bQ_{n-1}(t).$$

Multiplying by $-2\sqrt{a}$ yields

$$R''(t) + 4aR(t) = -bP'_{n-1}(t) + 2\sqrt{a}bQ_{n-1}(t).$$

Strictly speaking, this is again a harmonic equation, but we will not solve it as such. Since $R(t)$ is a polynomial, we shall look only for a polynomial solution. If $R(t) = \sum_{i=0}^{n-1} r_i t^i$, then

$$R''(t) = \sum_{i=0}^{n-3} r_{i+2}(i+2)(i+1)t^i.$$

It follows that

$$R''(t) + 4aR(t) = \sum_{i=0}^{n-3}(r_{i+2}(i+2)(i+1) + 4ar_i)t^i + 4ar_{n-2}t^{n-2} + 4ar_{n-1}t^{n-1}.$$

On the right-hand side we have the polynomial $-bP'_{n-1}(t) + 2\sqrt{a}bQ_{n-1}(t)$, equal to

$$\sum_{i=0}^{n-2}\left(-b\alpha[n-1,i+1](i+1) + 2\sqrt{a}b\beta[n-1,i]\right)t^i + 2\sqrt{a}b\beta[n-1,n-1]t^{n-1}.$$

To simplify the expressions, define

$$\gamma_i = -b\alpha[n-1,i+1](i+1) + 2\sqrt{a}b\beta[n-1,i], \quad 0 \le i \le n-2,$$

and

$$\gamma_{n-1} = 2\sqrt{a}b\beta[n-1,n-1].$$

Clearly, these coefficients will be redefined each time. Thus, we have $n$ numbers $\gamma_i, 0 \le i \le n-1$ and the system

$$\begin{cases} r_{i+2}(i+2)(i+1) + 4ar_i = \gamma_i, \ 0 \le i \le n-3 \\ 4ar_{n-2} = \gamma_{n-2} \\ 4ar_{n-1} = \gamma_{n-1} \end{cases}.$$

To solve, we also need the internal recurrence relation:

$$r_{n-1} = \tfrac{1}{4a}\gamma_{n-1}, \quad r_{n-2} = \tfrac{1}{4a}\gamma_{n-2},$$

and then, for all $r_k$ with $k$ ranging from $n-3$ down to 0, we obtain

$$r_k = \tfrac{1}{4a}(\gamma_k - r_{k+2}(k+2)(k+1)).$$

Thus the polynomial $R(t)$ is determined.

Recall that $S(t) = \tfrac{1}{2\sqrt{a}}(-bP_{n-1}(t) - R'(t))$. Therefore the coefficients of $S(t) = \sum_{i=0}^{n-1} s_i t^i$ are given by

$$s_i = \tfrac{1}{2\sqrt{a}}(-b\alpha[n-1,i] - (i+1)r_{i+1}), \quad 0 \le i \le n-2,$$

and
$$s_{n-1} = -\frac{b}{2\sqrt{a}}\alpha[n-1, n-1].$$

Now let us return to our problem, where we seek the polynomials $P_n(t), Q_n(t)$. We have essentially found them, up to the free coefficients. Indeed, since $P'_n(t) = R(t), Q'_n(t) = S(t)$, we obtain
$$\alpha[n, i] = \frac{r_{i-1}}{i}, \quad \beta[n, i] = \frac{s_{i-1}}{i}, \quad 1 \leq i \leq n.$$

It remains to compute $\alpha[n, 0], \beta[n, 0]$, which equal $P_n(0)$ and $Q_n(0)$, respectively.

Recall that $X_n(t) = P_n(t)\cos(\sqrt{a}t) + Q_n(t)\sin(\sqrt{a}t) + C_n$, where $C_n = \left(-\frac{b}{a}\right)^n C$. Moreover,
$$X'_n(t) = (P'_n(t) + Q_n(t)\sqrt{a})\cos(\sqrt{a}t) + (Q'_n(t) - P_n(t)\sqrt{a})\sin(\sqrt{a}t).$$

Since $X_n(0) = X_{n-1}(\tau)$ and $X'_n(0) = X'_{n-1}(\tau)$, we obtain
$$\begin{cases} P_n(0) + C_n = P_{n-1}(\tau)\cos(\sqrt{a}\tau) + Q_{n-1}(\tau)\sin(\sqrt{a}\tau) + C_{n-1}, \\ P'_n(0) + Q_n(0)\sqrt{a} = (P'_{n-1}(\tau) + Q_{n-1}(\tau)\sqrt{a})\cos(\sqrt{a}\tau) + (Q'_{n-1}(\tau) - P_{n-1}(\tau)\sqrt{a})\sin(\sqrt{a}\tau). \end{cases}$$

From the upper equation we find we find
$$\alpha[n, 0] = P_n(0) = P_{n-1}(\tau)\cos(\sqrt{a}\tau) + Q_{n-1}(\tau)\sin(\sqrt{a}\tau) + C_{n-1} - C_n.$$

In the lower equation we note that $P'_n(0) = \alpha[n, 1]$, already computed at this stage. Hence, since $\beta[n, 0] = Q_n(0)$, we obtain its value:
$$\frac{1}{\sqrt{a}}\left[(P'_{n-1}(\tau) + Q_{n-1}(\tau)\sqrt{a})\cos(\sqrt{a}\tau) + (Q'_{n-1}(\tau) - P_{n-1}(\tau)\sqrt{a})\sin(\sqrt{a}\tau) - P'_n(0)\right].$$

Thus the recurrence relation is obtained.

## 5.3 Explanations to Examples

The idea of our approach is to construct solutions $x_1(t), x_2(t)$ of delay equation (1.1) on corresponding interval and then to built equation (3.4). In the case of roots $\lambda_1$ and $\lambda_2$ satisfying inequalities $|\lambda_{1,2}| < 1$, we obtain the uniform exponential stability of delay equations. Construction of solutions can be done as we described.

Using recurrence relations, we can solve equations for every finite interval. We write a program which builds equation (3.4). Then it solves (3.4). Intervals, where $|\lambda_{1,2}| < 1$ give us the zones of the uniform exponential stability. Our program can be seen here: https://drive.google.com/drive/folders/10vyYZhaiy6OVugsmxn6jL-Dx90ByGG8e or using the QR-code:

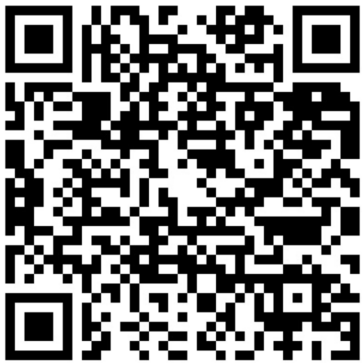

# 6 Discussion: conclusions and open problems

In this paper we propose a version of the Floquet theory for second order delay differential equations. Using standard definitions of the homogeneous equations, researchers came to infinite-dimensional fundamental systems. This does not allow to obtain full analog of the Floquet theory as it was noted in [31]. Various ideas to come to finite-dimensional version of the Floquet theory was proposed in many works [32, 34, 33, 47]. It can be noted that our version not only preserves the finite-dimensional fundamental system for delay equations, but also does not increase the order of the monodromy operators. Our approach develops the concept of act-and-wait control [48]. In our approach, existence of infinite number of intervals for a period $\omega$ such that equation is uniformly exponentially stable is obtained.

For future studies the following directions can be also considered. The case, where the period is less than the delay, was studied in the paper [33]. Its ideas in combination with ours can be a basis of a future research. The idea of our version of the Floquet theory can be developed for systems of delay equations using ideas of [48]. The recent paper [52], where act-and-wait control for delay systems was developed, could be a corresponding basis of a future research. In [52] the "act" period is divided into multiple segments, depending on the dimension of the open-loop system.

Using the idea of the Azbelev $W$-transform [4, 6, 8], we can obtain sufficient tests of stability for equations with sufficiently general variable delays. It looks, we could obtain results for equation with close to constant coefficients and delays. Note also one of possible directions for developments in the cases, where Assumptions 1 and 2 are not fulfilled. The way to achieve this can be constructing the Cauchy functions [4, 1, 8] for the corresponding equations with sampled data.

To compare our results with existing ones, we propose the following pictures for the equation, where the zone of stability obtained through D-division method [37] for the autonomous equations is red, and our method adds the blue domain of the uniform stability.

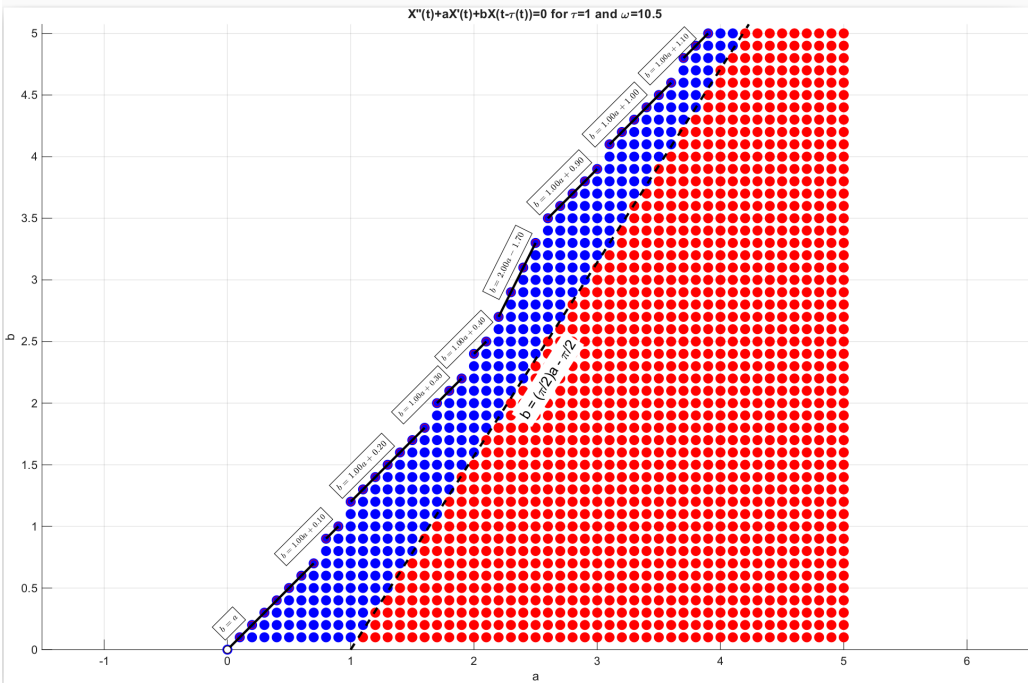

# 7 Acknowledgments

The work of the first author was partially supported by the grant "Oscillation Properties as a Key to Stability Studies", No. 0008174 of Science Forefront of the Ministry of Innovations, Science and Technology of the State of Israel. The paper is a part of PhD thesis of the third author.